\documentclass[11pt,reqno]{amsart}

\numberwithin{equation}{section}
\usepackage[all,cmtip]{xy}
\usepackage{amssymb,rotating,amsmath,amsfonts,amsthm,color,graphics,epsfig,bbm,bm,manfnt,psfrag}

\newcommand{\calL}{\mathcal{L}}

\newcommand{\calO}{\mathcal{O}}

\newcommand{\mB}{\mathbb{B}}
\newcommand{\mC}{\mathbb{C}}
\newcommand{\mD}{\mathbb{D}}

\newcommand{\mN}{\mathbb{N}}

\newcommand{\mR}{\mathbb{R}}

\newtheorem{theorem}{Theorem}[section]
\newtheorem{lemma}[theorem]{Lemma}

\theoremstyle{definition}

\newtheorem{remark}[theorem]{Remark}

\theoremstyle{definition}

\theoremstyle{definition}

\theoremstyle{definition}
\newtheorem{question}[theorem]{Question}

\theoremstyle{definition}
\newtheorem{conjecture}[theorem]{Conjecture}

\begin{document}



\title[Factorization in $SL_n(R)$]{Factorization in $SL_n(R)$ with elementary matrices when $R$ is the 
disk algebra and the Wiener algebra}

\author{Amol Sasane}
\address{Department of Mathematics, Faculty of Science, 
    Lund University, Sweden.}
\email{amol.sasane@math.lu.se}

\begin{abstract}
Let $R$ be the polydisc algebra or the Wiener algebra. It is shown that the group $SL_n(R)$ is 
generated by the subgroup of elementary matrices with all diagonal entries $1$ and 
at most one nonzero off-diagonal entry. The result an easy consequence of the deep result 
due to Ivarsson and Kutzschebauch~\cite{IvaKut}.
\end{abstract}

\maketitle

\section{Introduction}

Let $R$ be a commutative unital ring. Let $I_n$ denote the $n\times n$ identity matrix, 
that is the square matrix with all diagonal entries equal to $1\in R$ and off-diagonal entries equal to $0\in R$. 
Recall that an elementary matrix $E_{ij}(\alpha)$ over 
$R$ is a matrix of the form $I_n+\alpha \mathbf{e}_{ij}$, where $i\neq j$, $\alpha \in R$, and  
$\mathbf{e}_{ij}$ is the $n\times n$ matrix whose entry in the $i$th row and $j$th column is $1$ and all other entries are zeros. 
Let $SL_n(R)$ be the group of all $n\times n$ matrices whose entries are elements of $R$ and whose determinant is $1$. 
Let $E_n(R)$ be the subgroup of $SL_n(R)$ generated by the elementary matrices. 

A classical question in commutative algebra is the following: 

\begin{question}
\label{question}
Is $SL_n(R)$ equal to $E_n(R)$? 
\end{question}

The answer to this question depends on the ring $R$, and here is a list of a few known results.
\begin{enumerate}
\item If $R=\mC$, then the answer is ``Yes'', and this is standard exercise in linear algebra; 
 see for example \cite[Exercise~18.(c), page~71]{Art}.
\item Let $R$ be the polynomial ring $\mC[z_1, \cdots, z_n]$ in the indeterminates $z_1, \cdots, z_n$ 
with complex coefficients. 

\noindent If $n=1$, then the answer is ``Yes'', and this follows from the Euclidean Division 
Algorithm in $\mC[z]$. 

\noindent If $n=2$, then the answer is ``No'', and \cite{Coh} gave the following 
counterexample:
$$
\left[ \begin{array}{cc} 1+z_1 z_2 & z_1^2 \\
        -z_2^2 & 1-z_1 z_2 
       \end{array}\right] \in SL_2(\mC[z_1,z_2]) \setminus E_2(\mC[z_1,z_2]).
$$

\noindent For $n\geq 3$,  the answer is ``Yes'', and this is the $K_1$-analogue of Serre's Conjecture, 
which is the Suslin Stability Theorem \cite{Sus}.

\item The case of $R$ being a ring of continuous functions was considered in \cite{Vas}. Let $C(X; \mC)$ be 
the ring  of continuous complex-valued functions on 
the finite-dimensional normal topological space $X$ with pointwise operations. 
$C_b(X; \mC)$ denotes the subring of $C(X; \mC)$ consisting of {\em bounded} functions.  
It was shown in \cite{Vas} that for $R=C(X; \mC)$ or $C_b(X; \mC)$, the answer 
is ``Yes'' if there is no homotopy obstruction. Indeed, if $E$ is an elementary matrix, 
then $(X\owns) x\mapsto E(x) \in SL_n(\mC)$ is null-homotopic (to the constant map $x\mapsto I_n: X\rightarrow SL_n(\mC)$). 
 So it follows that  if $\pi(F)$ denotes the homotopy class of the map $x\mapsto F(x):X\rightarrow  SL_n(\mC)$ corresponding to $F\in SL_n(R)$, 
then a necessary condition for $F\in E_n(C(X; \mC))$ is that $\pi(F)=0$. It turns out that 
this condition is also sufficient, and this is the content of \cite[Theorem~4]{Vas}.

\item Based on the above result, it is natural to consider the question also 
 for the ring $\mathcal{O}(X)$ of holomorphic functions on Stein spaces in $\mC^n$. This was posed as an explicit open problem 
 by Gromov in \cite{Gro}, and was recently solved by Ivarsson and Kutzschebauch \cite{IvaKut}. The main result in 
 \cite{IvaKut} is the following:
 
 \begin{theorem}[\cite{IvaKut}]
 \label{theorem_IvaKut}
  If $X$ is a finite-dimensional reduced Stein space and $F:X\rightarrow SL_n(\mC)$ is a 
 holomorphic mapping that is null-homotopic, then there exists a natural number $K$ and holomorphic mappings 
 $G_1, \cdots, G_K: X\rightarrow \mC^{m(m-1)/2}$ such that $F$ can be written as a product of upper and lower 
 diagonal unipotent matrices 
 $$
 F(x)=M_1(G_1(x))\cdots M_K(G_K(x)), \quad x\in X, 
 $$
 where the matrices $M_j(G_j(x))$ are defined by 
 $$
 M_j(G_j(x)):= \left[ \begin{array}{ccc} 1 & & 0 \\  & \ddots & \\ G_j(x) & & 1\end{array}\right] 
 \textrm{ if } j\textrm{ is odd,}
 $$
 while 
 $$
 M_j(G_j(x)):= \left[ \begin{array}{ccc} 1 & & G_j(x) \\  & \ddots & \\ 0 & & 1\end{array}\right] 
 \textrm{ if } j\textrm{ is even.}
 $$
 In particular, the assumption of null-homotopy is always satisfied if $X$ is contractible. 
 \end{theorem}

\end{enumerate}

We wish to consider Question~\ref{question} for commutative, semisimple, unital complex Banach algebras $R$. 
 A special case is when $R=C_b(X; \mR)$, where $X$ is a compact Hausdorff topological space, and 
 item (3) above describes the answer in this special case. Motivated by this, we formulate the following question/conjecture, 
 but first we introduce some convenient notation. 
 
 Let $R$ be a commutative, semisimple, unital complex Banach algebra with maximal ideal space 
  denoted by $X_R$, equipped with the weak-$\ast$ topology induced from the dual space $R^*:=\calL(R; \mC)$ of $R$. 
  
  Let $\widehat{\cdot}: R \rightarrow C(X_R; \mC)$ denote the Gelfand transform. For $F\in SL_n(R)$, let $\widehat{F}$ be the 
  matrix with elements in $C(X_R; \mC)$ obtained by taking the 
  Gelfand transform of the entries of $F$, and $\pi(\widehat{F})$ denotes the homotopy class of 
   $ \varphi \mapsto \widehat{F}(\varphi): X_R \rightarrow SL_n(\mC).$

 \begin{conjecture}
 \label{conjecture}
  Let $R$ be a commutative, semisimple, unital complex Banach algebra.
  $F\in SL_n(R)$ belongs to $E_n(R)$ if and only if $\pi(\widehat{F})=0$.
 \end{conjecture}

We consider Question~\ref{question} for two important Banach algebras of holomorphic functions:
 the polydisc algebra $A(\overline{\mD}^n)$ and the Wiener algebra $W^+(\overline{\mD}^n)$. 

 Let $\mD:=\{z\in \mC: |z|< 1\}$ and $\overline{\mD}:=\{z\in \mC: |z|\leq 1\}$.
  Let $d\in \mN$. The {\em Wiener algebra}  $W^+(\overline{\mD}^n)$ is the Banach algebra defined by
 $$
 W^+(\overline{\mD}^d)\!=\!\left\{ \! \sum_{k_1=0}^\infty\! \cdots \!\sum_{k_d=0}^\infty \! a_{(k_1,\dots, k_d)} 
 z_1^{k_1}\cdots z_d^{k_d}:
 \sum_{k_1=0}^\infty\! \cdots\! \sum_{k_d=0}^\infty\! |a_{(k_1,\dots, k_d)}|<\infty\right\},
 $$
  with pointwise addition and multiplication, and the $\|\cdot\|_1$-norm given by
 $$
 \|f\|_{1}= \sum_{k_1=0}^\infty \!\cdots\! \sum_{k_d=0}^\infty\! |a_{(k_1,\dots, k_d)}|, \;\;
 f=\displaystyle \!\sum_{k_1=0}^\infty \!\cdots \!\sum_{k_d=0}^\infty \!a_{(k_1,\dots, k_d)} z_1^{k_1}\cdots z_n^{k_d}.
 $$
  The {\em polydisc algebra} $A(\overline{\mD}^d)$ is the Banach algebra
  of all continuous functions $f:\overline{\mD}^d\rightarrow \mC$ which are holomorphic in $\mD^d$,
  with pointwise addition and multiplication, and the supremum norm $\|\cdot\|_\infty$ given by
  $$
  \|f\|_\infty:=\sup_{(z_1,\dots, z_d)\in \mD^d}|f(z_1,\dots, z_d)|, \quad f\in A(\overline{\mD}^d).
  $$
 The ball algebra $A(\overline{\mB_d})$ is defined similarly, with the polydisc $\overline{\mD}^d$ replaced by the ball
  $$
 \overline{\mB_d}:=\{(z_1\cdots,z_d)\in \mC^d :|z_1|^2+\cdots+|z_d|^2\leq 1\}.
 $$
 For a $n\times n$ matrix $F$ with entries in $A(\overline{\mD}^d)$, $A(\overline{\mB_d})$ or $W^+(\overline{\mD}^d)$, we define 
 $$
 \|F\|:=\sum_{i,j =1}^n \|F_{ij}\|_{\infty}, 
 $$
 where $F_{ij}$ denotes the entry in the $i$th row and $j$th column of $F$. 
 Then $\|FG\| \leq \|F\| \|G\|$, for $n\times n$ matrices $F,G$ with entries from any of the 
 Banach algebras $A(\overline{\mD}^d)$, $A(\overline{\mB_d})$ or $W^+(\overline{\mD}^d)$.
 
 Our main result is the following.
 
 \begin{theorem}
 \label{main_theorem}
 If $R=A(\overline{\mD}^d)$, $A(\overline{\mB_d})$ or $W^+(\overline{\mD}^d)$, then 
  $SL_n(R)=E_n(R)$.  
 \end{theorem}
 
 If $R=A(\overline{\mD}^d)$ or $W^+(\overline{\mD}^d)$, then 
 in both cases, the maximal ideal space $X_R$ can be identified with $\overline{\mD}^d$ as a topological space. 
 Similarly, $X_{A(\overline{\mB_d})}=\overline{\mB_d}$. 
 If  Conjecture~\ref{conjecture} is true, then Theorem~\ref{main_theorem} follows from the 
 observation that $\overline{\mD}^d$, $\overline{\mB_d}$ are 
 contractible (since then $\pi(\widehat{F})$ is always trivial). 

 We will derive our main result as a consequence of the result from \cite{IvaKut} quoted above, and 
 \cite[Lemma~9]{Vas} reproduced below. 

\begin{lemma}[\cite{Vas}]
\label{lemma_9}
 Let $R$ be a commutative topological unital ring such that the set of 
 invertible elements of $R$ is open in $R$. If $F\in SL_n(R)$ is sufficiently close to $I_n$, 
 then $F$ belongs to $E_n(R)$. 
\end{lemma}

\section{Proof of Theorem~\ref{main_theorem}}
 
\begin{proof} We will simply prove the result in the case of the disc algebra $A(\overline{\mD}^d)$; the proofs 
in the cases of the ball algebra $A(\overline{\mB_d})$ and the Wiener algebra being analogous. 

Let $F\in SL_n(A(\overline{\mD}^d)$. Let $r\in(0,1)$ (to be determined later). Define 
$$
F_r(z_1, \cdots, z_d):=F(rz_1, \cdots , rz_d), \quad (z_1,\cdots, z_d)\in \mD^d.
$$
As $F_r \in \calO(\frac{1}{r}\mD^d)$, and $\det F_r\equiv 1$, it follows from Theorem~\ref{theorem_IvaKut} 
(since $\frac{1}{r}\mD^d$ is a contractible Stein domain) that there 
are elementary matrices $G_1, \cdots , G_K $ belonging to $E_n(\calO(\frac{1}{r}\mD^d))$ such that 
$$
F_r = E_1\cdots E_K \in E_n(\calO(\frac{1}{r}\mD^d)) \subset E_n(A(\overline{\mD}^d)). 
$$
Thus $F(I_n+ F^{-1} (F_r-F))=F_r \in E_n(A(\overline{\mD}^d))$. As $\det F=\det F_r=1$, it follows that also 
$\det (I_n+F^{-1} (F_r-F))=1$. We will be done if we manage to show that 
$I_n +F^{-1} (F_r-F) \in E_n(A(\overline{\mD}^d))$ too. But this is clear by Lemma~\ref{lemma_9}, since 
$$
\Big\|\Big(I_n +F^{-1} (F_r-F)\Big)- I_n\Big\|= \|F^{-1} (F_r-F)\|\leq \|F^{-1}\| \|F_r-F\|,
$$
and we can make $\|F_r-F\|$ as small as we like by choosing $r$ close enough to $1$. 
\end{proof}

\begin{remark}
 The above proof also works for some other Banach algebras of smooth functions 
 contained in the polydisc algebra, for example, if $N\in \mN$, the Banach algebra  
$\partial^{-N} A(\overline{\mD}^d)$ of all functions $f\in A(\overline{\mD}^d)$
whose complex partial derivatives of all orders up to $N$ belong to
$A(\overline{\mD}^d)$, with the norm
$$
\|f\|_{\scriptscriptstyle \partial^{-N}
  A(\overline{\mD}^d)}:=\sum_{\alpha_1+\dots+\alpha_d\leq N}\frac{1}{\alpha_1 ! \dots
  \alpha_d!}
\sup_{(z_1,\cdots, z_d)\in\mD^d}\left|\frac{\partial^{\alpha_1+\dots+\alpha_d}f}{\partial
    z_1^{\alpha_1} \dots \partial z_d^{\alpha_d}}(z_1,\cdots, z_d)\right|.
$$
\end{remark}

 In light of Theorem~\ref{main_theorem}, it is natural to ask the analogous question also for the Hardy algebra. 
 Recall that if $U$ is an open set in $\mC^d$, then the Hardy algebra $H^\infty(U)$ is the Banach algebra 
 of all complex-valued functions on $U$ that are bounded and holomorphic in $U$. 
 
\begin{conjecture}
 $SL_n(H^\infty(U))= E_n(H^\infty(U))$ if $U$ is the polydisc $\mD^d$ or open unit ball $U=\mB_d:=\{(z_1\cdots,z_d)\in \mC^d :
 |z_1|^2+\cdots+|z_d|^2<1\}$. 
\end{conjecture}

\end{document}